\documentclass{birkjour_t2}
\usepackage[utf8]{inputenc}
\usepackage{wrapfig}
\usepackage[font=small,labelfont=bf]{caption}
\newtheorem{thm}{Theorem}[section]

 \newtheorem{prop}[thm]{Proposition}
 \theoremstyle{definition}
 
 \theoremstyle{remark}

 \numberwithin{equation}{section}
\def\Pr{\noindent {\bf Proof. }}
\def\diver{\mathop{\mathrm{div}}\nolimits}
\def\tens#1{\mathbb{#1}}
\def\vec#1{\boldsymbol{#1}}
\def \R{\mathbb{R}}
\def \N{\mathbb{N}}
\def \O{\Omega}

\def\chiu{\hfill$\displaystyle\vspace{4pt}
\underset\Box\null$\par}

\def\bD{\tens{D}}

\def\bO{\tens{O}}
\def\bS{\tens{S}}

\def\bZ{\tens{Z}}

\def\bV{\tens{V}}

\def\bB{\tens{B}}
\def\bA{\tens{A}}
\def\bC{\tens{C}}

\def\bv{\vec{v}}
\def\bw{\vec{w}}
\def\bs{\vec{s}}

\def\bn{\vec{n}}
\def\bb{\vec{b}}
\def\bz{\vec{z}}
\def\bnul{\vec{0}}
\def\bphi{\vec{\varphi}}

\def\pa{\partial}
\def\na{\nabla}
\def\be{\begin{equation}}
\def\ba{\begin{array}}
\def\ea{\end{array}}
\def\ee{\end{equation}}
\def\displ{\displaystyle\vspace{6pt}}
\def\btau{\vec{\tau}}

\def\Lnd#1{L^{#1}_{\bn, \diver}}

\def\Wnd#1{W^{1,#1}_{\bn, \diver}}

\def\vs{\bv_\mathrm{s}}
\def\ff{_\mathrm{f}}
\def\rhotf{\varrho^\mathrm{m}\ff}
\def\rhots{\varrho^\mathrm{m}_\mathrm{s}}

\def\vf{\bv\ff}
\def\pa{\partial}

\def\Stot{\mathbb{S}}
\def\pf{p_\mathrm{f}}

\def\pft{p\ff^\mathrm{t}}

\begin{document}
\title{On unsteady flows of pore pressure-activated granular materials}
\author[A.~Abbatiello]{Anna Abbatiello}
\address{
Technische Universit\"{a}t Berlin, Institut f\"{u}r Mathematik,  Stra{\ss}e des 17. Juni 136,
10623 Berlin, 
 Germany} \email{anna.abbatiello@tu-berlin.de}
 \author[M.~Bul\'{i}\v{c}ek]{Miroslav Bul\'{i}\v{c}ek}
 \address{Charles University, Faculty of Mathematics and Physics, Mathematical Institute, Sokolovsk\'{a} 83, 18675 Prague~8, Czech Republic\\
 ORCiD: 0000-0003-2380-3458} \email{mbul8060@karlin.mff.cuni.cz}
\author{Tom\'{a}\v{s} Los}
 \address{Charles University, Faculty of Mathematics and Physics, Mathematical Institute, Sokolovsk\'{a} 83, 18675 Prague~8, Czech Republic}
\email{los@karlin.mff.cuni.cz}
\author{Josef M\'{a}lek}
 \address{Charles University, Faculty of Mathematics and Physics, Mathematical Institute, Sokolovsk\'{a} 83, 18675 Prague~8, Czech Republic\\
 ORCiD: 0000-0001-6920-0842}
\email{malek@karlin.mff.cuni.cz}
\author{Ond\v{r}ej~Sou\v{c}ek}
 \address{Charles University, Faculty of Mathematics and Physics, Mathematical Institute, Sokolovsk\'{a} 83, 18675 Prague~8, Czech Republic}
\email{soucek@karel.troja.mff.cuni.cz}
\date{}
\begin{abstract}
We investigate mathematical properties of the system of nonlinear partial differential equations that describe, under certain simplifying assumptions, evolutionary processes in water-saturated granular materials. The unconsolidated solid matrix behaves as an ideal plastic material before the activation takes place and then it starts to flow as a Newtonian or a generalized Newtonian fluid. The plastic yield stress is non-constant and depends on the difference between the given lithostatic pressure and the pressure of the fluid in a pore space. We study unsteady three-dimensional flows in an impermeable container, subject to stick-slip boundary conditions. Under realistic assumptions on the data, we establish long-time and large-data existence theory.
\end{abstract}
\maketitle
\section{Introduction} 
The purpose of this study is to investigate mathematical properties of a system of nonlinear partial differential equations (PDEs) developed in \cite{old} to describe processes, such as static liquefaction or enhanced oil recovery, in water-saturated (geological) materials. Such materials can be viewed as two component mixtures consisting of a granular solid matrix and a fluid filling the interstitial pore space. More specifically, we investigate the following system of PDEs:
\begin{subequations}
\label{fullsystem}
\begin{align}
\label{eq1}
\diver\bv &= 0,\\
\label{}
\rhots\left({\pa_t\bv} + \diver{(\bv{\otimes}\bv)}\right) &= \diver\Stot -\nabla p + \rhots\bb\ ,\\
{\pa_t\pf} + \bv\cdot\nabla\pf &= K\Delta\pf  - \diver(K\rhotf\bb) + \pa_t p_s + \bv\cdot\nabla p_s\ ,\label{eq3}\\
 \vf &= \bv - \frac{1}{\alpha}\widehat{\phi}(p-\pf)\left(\nabla\pf - \rhotf\bb \right),\label{vf}
\end{align}
where $\bS$ and $\bD\bv$ satisfy
\begin{equation} \label{S2}
\begin{split}
 \bD\bv=\bO &\Rightarrow \ |\bS|\leq \tau(\pf),\\
\bD\bv\neq\bO &\Rightarrow \  \bS =\tau(\pf)\frac{ \bD\bv}{|\bD\bv|}+2 \nu_* \left(|\bD\bv|-\delta_*\right)^+\frac{\bD\bv}{|\bD\bv|}, 
\end{split}
\hspace{6pt}\mbox{   with }\tau({\pf}) := q_*(p_s- \pf)^+.\end{equation} 
\end{subequations}
The system \eqref{fullsystem} coincides with equations (2.24)--(2.26) stated in \cite{old} (see also \cite{CM}) provided that we set $\delta_*=0$ in \eqref{S2} and we identify the symbols $\bv$ and $\pf$ with $\vs$ and $\pft$ used in \cite{old}. In  \cite{old}, equations \eqref{fullsystem} are summarized at the end of Section 2 as the outcome of derivation starting from the general principles of the theory of interacting continua, also using several well-motivated simplifying assumptions. In \eqref{fullsystem}, $\bv$ represents the velocity of  the granular solid matrix, $\vf$ is the velocity of the interstitial fluid, $p$ stands for the total pressure of the whole mixture and $\pf$ is the pressure of the fluid in a pore space. The vector and scalar fields $\bv, \vf, p$ and $\pf$ represent the unknowns, the other quantities are given material functions/parameters. More precisely, $\phi=\widehat{\phi}(p-\pf)$ is the porosity given as a function of the 
``effective'' pressure $p-\pf$, $\rhots$ and $\rhotf$ are the constant material densities of the solid and the fluid, $\bb$ represents given external forces, $\alpha$ is the drag coefficient, $\nu_*$ the viscosity of the fluid, $K$ is a constant coefficient related to permeability, $\delta_*$ is the non-zero activation parameter and $p_s$ is the given lithostatic pressure. Note that $\vf$ appears only in \eqref{vf} and can be always obtained a posteriori once $\bv, \pf$ and $p$ are obtained from \eqref{eq1}--\eqref{eq3}. Consequently, in what follows, we consider system \eqref{fullsystem} without equation \eqref{vf}. It is worth observing that the constitutive relation \eqref{S2} can be rewritten in a more compact way as an implicit constitutive relation (see Figure 1):
\begin{equation}\label{eqS}
\bS= \bZ+2 \nu_* \left(|\bD\bv|-\delta_*\right)^+\frac{\bD\bv}{|\bD\bv|} \mbox{ with } \bZ \mbox{ fulfilling } (|\bZ|\!-\!\tau(\pf))^+ \!\!+\! ||\bD\bv|\bZ\!-\!\tau(\pf)\bD\bv|\!=\!0. 
\end{equation}
We will exploit formulation \eqref{eqS} in our analysis. A systematic study of implicit constitutive equations go back to the original works \cite{RA} and \cite{RA2}.

We study the system of PDEs \eqref{eq1}--\eqref{eq3} and \eqref{eqS} in time-space cylinder $(0, T)\times\Omega$, where $T>0$ and $\Omega\subset \mathbb{R}^3$ is a bounded open connected set with Lipschitz boundary $\partial\Omega$.
 We  complete the system by considering the following boundary and initial conditions:   
\begin{subequations}\label{ibc}\begin{align}
&\!\bv\cdot\bn=0 \mbox{ and }  \na \pf \cdot \bn=0  \mbox{ on } (0,T)\times\partial\Omega, \label{boundary0} \vspace{6pt}\\ \vspace{6pt}
&\!\bs= \bz+\gamma_* \left(|\bv_\tau|-\beta_*\right)^+ \frac{\bv_\tau}{|\bv_\tau|} \mbox{ with } \bz \mbox{ fulfilling }(|\bz|\!-\!s_*)^+ \!+\! ||\bv_\tau|\bz-s_*\bv_\tau|\!=\!0  \mbox{ on } (0,T)\times\partial\Omega, 
\label{boundary2}\\ \vspace{6pt}
&\!\bv(0, \cdot)=\bv_0  \mbox{ and }    \pf(0,\cdot)={p}_0  \mbox{ in }  \O. \label{initial}
\end{align}\end{subequations}
Here, we used the following notation: $\bn :\pa\O\to \mathbb{R}^3$ stands for the unit outer normal vector, 
\begin{wrapfigure}{h!}{0.28\textwidth}
  \vspace{-6pt}
\includegraphics[scale=0.9]{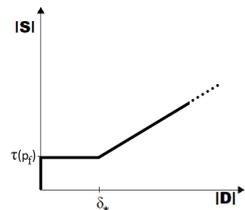}
  \vspace{-6pt}
\caption[width=5cm]{\footnotesize{Representation of the material response described by \eqref{eqS}}}
\label{fig:fig1}
\vspace{-10pt}
\end{wrapfigure}
while for any vector $\bz$ defined on $\pa\O$,  $\bz_\tau:=\bz-(\bz\cdot \bn)\bn$ denotes the tangential component of $\bz$, 
in particular, $\bs:=-(\bS\bn)_\tau$, and  $\gamma_*, \beta_*, s_*$ are non-negative constants. Condition \eqref{boundary2} describes
the shifted stick-slip (or threshold slip) and it is analogous to that for the stress tensor in the bulk (see \eqref{eqS}). It includes as special cases, the stick-slip  by taking $\beta_*=0$ while  $s_*, \gamma_*>0$, Navier's slip $\bs=\gamma_*\bv_\tau$ by taking $s_*, \beta_*=0$ while $\gamma_*>0$, and perfect slip $\bs=\bnul$ by setting $s_*=\gamma_*=0$. Note that the no-slip condition is obtained by letting either $s_*\to +\infty$ or by setting $\beta_*=0$ and letting $\gamma_*\to +\infty$.

The main purpose of this study is to establish long-time and large-data theory to the initial- and boundary-value problem described by \eqref{eq1}--\eqref{eq3},  \eqref{eqS}, \eqref{ibc}, see Theorem \ref{main-thm1} below. The novelties consist not only in incorporating a more general model with $\delta_*\geq0$, but more importantly in providing a different proof for more general class of data (particularly for $\bb$ that is merely $L^2$-integrable). More precisely, we can avoid using $L^\infty$-estimates for $\pf$ needed in~\cite{old}. Consequently, the main tool for taking the limit in the constitutive equations cannot be applied in the form given in \cite[Proposition 5.3]{old}, but has to be modified in an essential way due to a lower integrability of $\pf$, but also a more complicated material response. 
\begin{wrapfigure}{h!}{0.28\textwidth}
  \vspace{-6pt}
\includegraphics[scale=0.3]{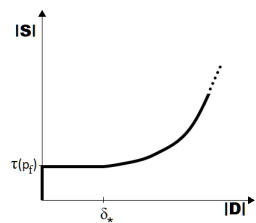}
\vspace{-10pt}
\caption{\footnotesize{Representation of the material response described by \eqref{S2new}, whenever $q>2$}}
\label{fig:fig2}
\vspace{-6pt}
\end{wrapfigure}

The novel key tool regarding the attainment of the constitutive equations by the limiting objects is proved separately in Proposition \ref{convergence-lemma}. The key assumption of this proposition, namely \eqref{limsupZ} and \eqref{limsupV}, call for taking $\bv^n-\bv$ as a test
function in the weak formulation of balance of linear momentum. However, $\bv^n-\bv$ is not admissible test function in the setting considered here. 
This difficulty can be overcome  by using  the $L^\infty$-truncation method which requires to introduce an integrable pressure, as the truncations $(\bv^n-\bv)_\infty$ are not divergenceless. 
Following the approach originally developed in \cite{BMR} (see also \cite{BM}), we
overcome such difficulty by considering slipping boundary conditions \eqref{boundary0}--\eqref{boundary2}. 
As pointed out in \cite{BleMalRaj}, the analysis for unsteady flows changes remarkably when the no--slip condition is considered.   

We use the $L^\infty$-truncation method in proving Theorem \ref{main-thm1} below. While the truncations
 $(\bv^n-\bv)_\infty$ are difficult to make solenoidal, the authors of \cite{sol-lip} succeeded to make the Lipschitz approximations $(\bv^n-\bv)_{1, \infty}$ divergenceless and they thus developed a solenoidal version of the Lipschitz truncation method. This tool allows one to avoid the presence of the pressure in the setting, therefore one may include more general responses as well as boundary conditions. 
As a matter of fact, we  present  new results available for systems describing materials that behave after activation $|\mathbb{D}\bv|>\delta_*$, as a power-law fluid, i.e. the constitutive equation \eqref{eqS} is replaced by (see Figure 2)
\begin{equation}\label{S2new}
\bS= \bZ+ 2 \nu_*|\bD\bv|^{q-2}\!\left(|\bD\bv|-\delta_*\right)^+\frac{\bD\bv}{|\bD\bv|}  \mbox{ with } \bZ \mbox{ fulfilling }(|\bZ|-\tau(\pf))^+ \!\!+\! ||\bD\bv|\bZ-\tau(\pf)\bD\bv|\!=\!0.
\end{equation}
The available results are presented in Theorem \ref{main-thm2}.  We are not providing the proof of these results as they can be deduced from the approach used when proving Theorem \ref{main-thm1} and from the methods used recently, for example, in \cite{BleMalRaj}. Note that the latter results are restricted to models \eqref{S2new} with $q>\frac{6}{5}$ (in three dimensions). Recently, another concept of dissipative solution was introduced in \cite{AbbFeir} and, its long-time and large-data existence is proved independently of what is the value of $q$ (in particular also for $q\in[1, {6}/{5}]$). 
In fact in the theory developed in \cite{AbbFeir} the stress tensor can be merely subdifferential of a convex potential depending on $\mathbb{D}\bv$, whose growth is at least linear.  There are other approaches to analyze the mathematical properties of  Bingham fluids (see e.g. \cite{MaRuSh} and  \cite{She}), but they are usually based on regularity techniques requiring smoother data.

\section{Preliminaries and main results}
For the sake of simplicity in the  right-hand side of \eqref{eq3}, which  has the form $g:=\pa_t p_s -\diver \bb$, we omit the effect of $\pa_t p_s$ as it plays the role of a given external force and it can be easily incorporated into the analysis. We also set without loss of any generality $\rhots=\rhotf=K= 2\nu_* = \gamma_* =q_*=1,$ while we assume $\delta_*, s_*, \beta_*\geq0$. Finally, to shorten the notation we set  $Q:=(0,T)\times \O$ and $\Sigma:= (0,T)\times\pa\O$, where we fix $\O$ to be a bounded open set in $\mathbb{R}^3$ with either Lipschitz or $C^{1,1}$ boundary $\partial \Omega$; such sets are denoted either $\Omega \in C^{0,1}$ or $\Omega \in C^{1,1}$.

Before stating the main results, let us summarize the notation. The symbol $\bD\boldsymbol\varphi$ stands for the symmetric part of the gradient of a vector-valued function $\boldsymbol\varphi$, i.e. 
$$\bD\boldsymbol\varphi:=\frac{\nabla\boldsymbol\varphi+(\nabla\boldsymbol\varphi)^T}{2}.$$
The symbols $(L^q(\O), \|\cdot\|_q)$ and $(W^{1,q}(\O), \|\cdot\|_{1,q})$ with $q\in [1,\infty]$, stand respectively for the Lebesgue spaces, the  Sobolev spaces with their own norms.
If $X$ is a Banach space of scalar functions, then $X^3$ and $X^{3\times3}$ denote the space of vector-valued functions having three components and  the space of tensor-valued functions  respectively, with each component belonging to $X$. For a Banach space $X$,   $L^q(0,T; X)$ denotes a corresponding Bochner space. We make use of the following function spaces
\begin{align*}\displ
&\Lnd{q} :=  \overline{\left\{ \bv \in C^\infty(\Omega)^3; \, \diver \bv =0 \ \mbox{in} \ \O; \, \bv\cdot \bn=0 \ \mbox{on} \ \pa \O\right\}}^{\|\cdot\|_{L^q(\Omega)^3}} \mbox{ for } q\in [1,\infty),\\
\displ
&W_{\bn}^{1,q}  :=  \{\bv\in W^{1,q} (\Omega)^3; \bv\cdot \bn = 0 \ \mbox{on}\ \partial \Omega\},\\\displ
&W_{\bn,\diver}^{1,q}  :=  \{\bv\in W^{1,q} (\Omega)^3; \, \diver \bv=0\ \mbox{in}\ \Omega; \, \bv\cdot \bn=0 \ \mbox{on} \ \partial \Omega\},
\end{align*}
while ${(W_{\bn}^{1,q})}^*$, ${(W_{\bn,\diver}^{1,q})}^*$ are the dual spaces to $W_{\bn}^{1,q}$  and $W_{\bn,\diver}^{1,q}$ respectively.
In particular, assuming  $\O\in C^{1,1}$  the following Helmholtz decomposition holds
$$
W^{1,2}_{\bn}=W^{1,2}_{\bn,\diver} \oplus \{\nabla \varphi; \varphi \in W^{2, 2}(\Omega), \nabla \varphi \cdot \bn =0 \textrm { on } \partial \Omega\}.
$$
Note that such decomposition is not valid for $W^{1,2}_{0}(\Omega)^3$.  
\\

We are ready to enunciate the first result, which is the existence of weak solutions to system \eqref{eq1}--\eqref{eq3},  \eqref{eqS}, \eqref{ibc}, proved in Section \ref{sectionProof}.

\begin{thm}\label{main-thm1}
For any $\O\in C^{1,1}$, $T > 0$ and for any  $ \bv_0, p{_0}, \bb, p_s$ fulfilling 
$$\bv_0\in \Lnd{2}, \ p{_0}\in L^2(\O), \ \bb\in L^2(Q), \ p_s\in L^5(Q),$$
there exists a quintuplet $(\bv,\pf, p, \bS, \bs)$:
 \begin{align*}
&\bv \!\in\!  L^\infty(0,T; \Lnd{2})\!\cap\! L^2(0,T;W^{1,2}_{{\bf n},{\rm div}}),   \partial_t\bv\!\in\! {(L^2(0,T;W^{1,2}_{\bn})\!\cap \!L^5(Q)^3)}^*,\\
&\pf \in L^\infty(0,T; L^2(\O))\cap L^2(0,T; W^{1,2}(\O)), \ \partial_t \pf \in (L^4(0,T; W^{1,2}(\O)))^*,\\
& p=p_1+p_2 \mbox{ where } p_1\in L^2(Q) \mbox{ and } p_2 \in L^{\frac{5}{4}}(0,T; W^{1, \frac{5}{4}}(\O)), \\
&\bS\in  L^{2}(Q)^{3\times 3}, \ \bs\in  L^{\frac{8}{3}}(\Sigma)^3,
\end{align*}
satisfying the following weak formulations: 

\begin{align}\!\!\!\!\!\begin{array}{l}\displ\vspace{6pt}
\int_0^T\!\!\langle \partial_t\bv, \bw\rangle + \int_{Q}\bS:\!\bD\bw
 -\int_{Q} (\bv\otimes \bv)\!:\! \bD\bw +\int_{\Sigma} \bs\cdot\bw_{\btau} = \int_{Q}\bb\cdot \bw
  + \int_{Q}\! p_1\diver \bw - \int_{Q}\! \nabla p_2\cdot \bw 
   \\\displ\vspace{6pt}
\hfill \mbox{ for all } \bw\in L^2(0,T;W^{1,2}_{\bn})\!\cap \!L^5(Q)^3,
\end{array}\\ \vspace{6pt}
\int_0^T\!\!\! \langle \partial_t \pf, z\rangle -\!\! \int_{Q}\pf \bv\cdot \nabla z + \int_{Q}\nabla \pf\cdot\!\nabla z =\!\! \int_{Q}\bb\cdot\nabla z  - \int_{Q} p_s\bv\cdot \nabla z \
  \mbox{ for all } z\!\in\! L^4(0,T; W^{1,2}(\O)),
\end{align}
and the following constitutive equations:

\begin{align}
&\!\!\!\!\!\!\!\begin{array}{l}\displ
\bS =  \bZ\!+\! \left(|\bD\bv|-\delta_*\right)^+\!\frac{\bD\bv}{|\bD\bv|}  \mbox{ with } \bZ \mbox{ fulfilling }(|\bZ|-\tau(\pf))^+\! + ||\bD\bv|\bZ-\tau(\pf)\bD\bv|\!=\!0\\ \displ\vspace{6pt}
\hfill  \mbox{ with } \tau({\pf})\! =\! (p_s- \pf)^+ \mbox{ a.e. in } Q,\end{array}
\\ \label{obiettivo-boundary}
&\!\!\!\!\!\bs= \bz+ \left(|\bv_\tau|-\beta_*\right)^+ \frac{\bv_\tau}{|\bv_\tau|} \mbox{ with } \bz \mbox{ fulfilling } 
(|\bz|-s_*)^+ + ||\bv_\tau|\bz-s_*\bv_\tau|=0  \mbox{ a.e. on } \Sigma,
\end{align}
and attaining the initial conditions in the following sense: 

\begin{equation}\label{initial-weak} \lim_{t\rightarrow 0+}\|\bv(t)-\bv_0\|_2=0, \  \  \lim_{t\rightarrow 0+}\|\pf(t)-p_0\|_2=0.
\end{equation}
\end{thm}

The second result concerns system  \eqref{eq1}--\eqref{eq3}, \eqref{ibc}, and \eqref{S2new}.
\begin{thm}\label{main-thm2}
Let $\O \in C^{0,1}$, $T > 0$, and $q>\frac{6}{5}$. Set $m:=\max\{2, q'\}$ and $r:=\max\left\{ q, \frac{5q}{5q-6}\right\}.$ For any  $ \bv_0, p{_0}, \bb, p_s$ fulfilling 
$$\bv_0\in \Lnd{2}, \ p{_0}\in L^2(\O), \ \bb\in L^{m}(Q), \ p_s\in L^{\frac{10q}{5q-6}}(Q),$$
 there exists a quadruplet $(\bv,\pf, \bS, \bs)$:
 \begin{align*}
&\bv \in  L^\infty(0,T; \Lnd{2}) \cap L^q(0,T;W^{1,q}_{{\bf n},{\rm div}}),  \partial_t\bv\in L^{r'}(0, T; (\Wnd{r})^*),\\
&\pf \!\in\! L^\infty(0,T; L^2(\O))\!\cap\! L^2(0,T; W^{1,2}(\O)), \partial_t \pf\! \in\!(L^2(0,T; W^{1, 2}(\Omega))\!\cap\! L^{\frac{10q}{5q-6}}(Q))^*,\\
&\bS\in  L^{q'}(Q)^{3\times 3},  \ \bs\in  L^{2}(\Sigma)^3,
\end{align*}

satisfying the initial conditions \eqref{initial-weak}, the following weak formulations: 

\begin{equation}
\!\! \int_0^T\!\!\!\langle \partial_t\bv, \bphi\rangle+
 \int_{Q}\!\!\!\bS: \bD\bphi
 - \int_{Q}\!\!(\bv\otimes \bv)\!:\! \bD\bphi
  +\int_{\Sigma} \bs\cdot\bw_{\btau}=\! \int_{Q}\!\!\bb\cdot \bphi
\ \ \ \textrm{ for all } \bphi\in L^r(0,T; \Wnd{r}), 
\end{equation}

\begin{equation}
\!\!\! \int_0^T\!\!\!\langle \partial_t\pf, \varphi\rangle +\int_{Q}\!\!\nabla \pf\cdot \nabla \varphi - \int_{Q}\! \pf \bv\cdot \nabla\varphi =\!\!\int_{Q}\!\!(\bb-p_s\bv)\cdot\nabla \varphi \ \ \
\mbox{ for all }  \varphi\! \in\! L^2(0,T; W^{1, 2}(\Omega))\cap L^{\frac{10q}{5q-6}}(Q),
\end{equation}

and the following constitutive equations:
\begin{align}
&\!\!\!\!\!\!\!\!\!\!\!\begin{array}{l}\displ\bS= \bZ+ \left(|\bD\bv|-\delta_*\right)^+|\bD\bv|^{q-2}\frac{\bD\bv}{|\bD\bv|}  \mbox{ with } \bZ \mbox{ fulfilling }
(|\bZ|-\tau(\pf))^+ + ||\bD\bv|\bZ-\tau(\pf)\bD\bv|=0\\ \displ\vspace{6pt}\hfill \mbox{ with } \tau({\pf}) = (p_s- \pf)^+ \textrm{  a.e. in } Q,
\end{array}\\
&\!\!\!\!\!\!\!\!\bs= \bz+ \left(|\bv_\tau|-\beta_*\right)^+ \frac{\bv_\tau}{|\bv_\tau|} \mbox{ with } \bz \mbox{ fulfilling } 
(|\bz|-s_*)^+ + ||\bv_\tau|\bz-s_*\bv_\tau|=0  \mbox{ a.e. on } \Sigma. 
\end{align}
\end{thm}

This result is stated without the proof here. The proof can be however achieved in the spirit of Theorem \ref{main-thm1} by  employing a solenoidal version of the Lipschitz-truncation method developed in  \cite{sol-lip}, and by using the approximation scheme presented in \cite[Theorem 3.3]{BleMalRaj}. 

\section{Attainment of the constitutive equations}
In this section, we establish a new scheme how to take the limit in the constitutive equations needed when proving Theorem \ref{main-thm1}.

\begin{prop}\label{convergence-lemma}
Let $U\subset Q$ be an arbitrary measurable bounded set and let $\{\bZ^n\}_{n=1}^{+\infty}$,  
$\{\bD^n\}_{n=1}^{+\infty}$ and $\{\pf^n\}_{n=1}^{+\infty}$
be sequences such that 
\begin{align}
&\bZ^n=\tau(\pf^n) \frac{\bD^n}{|\bD^n|+\frac{1}{n}} \mbox{ with } \tau(\pf^n)= (p_s-\pf^n)^+ \mbox{ a.e. in }U, \label{An}\\ 
&\bZ^n \rightharpoonup \bZ \ \textrm{weakly in } L^{2}(U)^{3\times 3}, 
 \label{Zw}\\\displ
&\bD^n \rightharpoonup \bD \ \textrm{weakly in } L^{2}(U)^{3\times 3}, \label{Dw} \\\displ
&\pf^n  \to \pf  \textrm{ strongly in }  L^2(U) \mbox{ and a.e. in }U,\label{ps}\\\displ
&\limsup_{n\to \infty} \int_{U} \bZ^n : \bD^n \leq \int_{U} \bZ : \bD. \label{limsupZ}
\end{align}
Then
\begin{equation}\vspace{6pt}
 (|\bZ|\!-\!\tau(\pf))^+ \!\!+\! ||\bD|\bZ\!-\!\tau(\pf)\bD|\!=\!0.  \label{C1} 
\end{equation}
In addition, assume that $\{\bV^n\}_{n=1}^{+\infty}$ is a sequence such that
\begin{equation}
\bV^n=\left(1-\frac{\delta_*}{|\bD^n|}\right)^+\bD^n \ \ \mbox{ a.e. in }U, \label{tildeAn}
\end{equation}
fulfilling
\begin{align}\displ \vspace{6pt}
&\bV^n \rightharpoonup \bV \ \textrm{weakly in } L^{2}(U)^{3\times 3},\vspace{6pt}
 \label{Vw}\\\displ \vspace{6pt}
 &\limsup_{n\to \infty} \int_{U} \bV^n : \bD^n \le \int_{U} \bV : \bD.  \label{limsupV}\vspace{6pt}
 \end{align}
 Then 
\begin{equation}
\bV= \left(1-\frac{\delta_*}{|\bD|}\right)^+\bD \ \ \mbox{ a.e. in }U.\label{C3}
\end{equation}
 \end{prop}
 
\Pr 
First, note that by  virtue of \eqref{ps} and the Lipschitz-continuity of $\tau$, it follows that 
\begin{equation}\label{tau}
\tau(\pf^n) \to \tau(\pf) \mbox{ strongly in } L^2(U).
\end{equation}

Now, as for any  $\bA\in L^2 (U)$, $\bA\neq\mathbb{O}$,  it holds
\begin{equation}\label{mon}
\left(\tau(\pf) \frac{\bD^n}{|\bD^n|+ \frac{1}{n}} - \tau(\pf) \frac{\bA}{|\bA|+\frac{1}{n}}\right): \left(\bD^n -\bA\right)\geq 0,
\end{equation}
integrating this inequality  over $U$,  subtracting and adding $\bZ^n$ and using  \eqref{An}, 
we get 

\begin{equation}\label{J1J2}
\int_{U} \left(\tau(\pf) \frac{\bD^n}{|\bD^n|+\frac{1}{n}} -\tau(\pf^n) \frac{\bD^n}{|\bD^n|+\frac{1}{n}}\right): \left(\bD^n -\bA\right) \\
+ \int_{U}\left( \bZ^n - \tau(\pf) \frac{\bA}{|\bA|+\frac{1}{n}}\right):\left(\bD^n - \bA\right)\geq0.
\end{equation}

Taking limsup as $n\to \infty$ and employing the facts that the first integral converges to zero due to \eqref{tau}, the term $\frac{\bD^n}{|\bD^n|+\frac{1}{n}}$ is uniformly bounded in $L^\infty (U)$ and  the sequence  $\bD^n -\bA$  is bounded in $L^2 (U)$,  we obtain
\begin{equation}\label{J2}
 \limsup_{n\to \infty}\int_{U} \left( \bZ^n - \tau(\pf) \frac{\bA}{|\bA|+\frac{1}{n}}\right):\left(\bD^n - \bA\right)\geq0.
 \end{equation}

Referring then to the convergences \eqref{Zw} and \eqref{Dw} and using also \eqref{limsupZ}, we conclude that 

\begin{equation}\label{preMinty} 
\int_{U} \left( \bZ - \tau(\pf) \frac{\bA}{|\bA|}\right):\left(\bD - \bA\right)\geq0.
\end{equation}

Now, for any $\delta>0$, $\varepsilon \in (0, \delta)$ and for arbitrary  matrices $\bC$ and $\bB_1$ bounded in $L^2(U)^{3\times 3}$ and satisfying $|\bC|\le 1$ and $ \bB_1 \neq \bO$,  consider   
$$
\bA:= \bB_1 \, \chi_{\{|\bD|=0\}} + (\bD - \varepsilon \bC) \, \chi_{\{|\bD|>\delta\}} + \bD \, \chi_{\{0<|\bD|\le \delta\}}. 
$$  
Note that such $\bA$'s are non-zero in $U$. Inserting them into \eqref{preMinty} we obtain 
\begin{equation}\label{preMinty_p} 
- \int_{\{|\bD|=0\}} \left( \bZ - \tau(\pf) \frac{\bB_1}{|\bB_1|}\right):\bB_1 + \varepsilon \int_{\{|\bD|>\delta\}} \bC: \left( \bZ - \tau(\pf) \frac{\bD- \varepsilon \bC}{|\bD- \varepsilon \bC|}\right) \geq 0.
\end{equation}
Letting first $\varepsilon \to 0$ in \eqref{preMinty_p}, we observe that 
\begin{equation}\label{abba}
  \int_{\{|\bD|=0\}} \bZ : \bB_1 \le \int_{\{|\bD|=0\}} \tau(p_f) |\bB_1| 
\end{equation}
for any $\bB_1\neq \mathbb{O}$. Consider, for any $a>0$ and $\omega\subset U$, the matrix $\bB_1$ of the form 
$$\bB_1=a\, \mathbb{I}\chi_{\{(U\setminus\omega) \cup\{\bZ=0\}\}}+ \frac{\bZ}{|\bZ|}\chi_{\{\omega\setminus\{\bZ=0\}\}}.$$
It then follows from \eqref{preMinty_p} that
$$ \int_{\{|\bD|=0\}\cap\, \omega\,\cap \{\bZ\neq 0\}}\!\!\! |\bZ|  \le \int_{\{|\bD|=0\}\cap\, \omega\cap \{\bZ\neq 0\}}\!\!\! \tau(\pf) + a\, C\!\int_{(U\setminus\omega) \,\cup\{\bZ=0\}}\!\!(\tau(\pf)+ |\bZ|)$$
with $C$ positive constant, which implies letting $a\to 0$
$$\int_{\{|\bD|=0\}\cap\, \omega\,\cap \{\bZ\neq 0\}} |\bZ|  \le \int_{\{|\bD|=0\}\cap\, \omega\cap \{\bZ\neq 0\}} \tau(\pf).$$
Since $\omega$ is arbitrary, we conclude that
\begin{equation}\label{abba2}
|\bZ|\leq \tau(\pf) \mbox{ on the set } \{|\bD|=0\}.
\end{equation}

Next, letting $|\bB_1|\to 0$ in \eqref{preMinty_p}, employing \eqref{abba}, we get 
\begin{equation*} 
\int_{\{|\bD|>\delta\}} \bC: \left( \bZ - \tau(\pf) \frac{\bD- \varepsilon \bC}{|\bD- \varepsilon \bC|}\right) \geq 0,
\end{equation*}
which, after letting $\varepsilon \to 0$, leads to 
\begin{equation*} 
 \int_{\{|\bD|>\delta\}} \bC: \left( \bZ - \tau(\pf) \frac{\bD}{|\bD|}\right) \geq 0.
\end{equation*}
Finally, letting $\delta\to 0$, we get, for arbitrary $\bC$,  
$$ \int_{\{|\bD|>0\}} \bC: \left(\bZ -\tau(\pf) \frac{\bD}{|\bD|}\right)\geq 0.$$ 
This implies 
\begin{equation}\label{abba3}
\bZ =\tau(\pf) \frac{\bD}{|\bD|} \mbox{ when } |\bD|\neq 0.
\end{equation}
The latter and \eqref{abba2} are equivalent to \eqref{C1}.
\\

It remains to prove \eqref{C3}, which however follows from standard Minty's argument. Indeed, by the monotonicity, we have 
$$\limsup_{n\to \infty}\int_{U}\left(\bV^n-\bA\left(1-\frac{\delta_*}{|\bA|}\right)^+\right):\left(\bD^n-\bA\right)\geq0$$
for any $\bA\in L^2(U)^{3\times 3}$. By virtue of \eqref{limsupV} and of  convergences \eqref{Vw} and \eqref{Dw} we get 
 
$$ \int_{U}\left(\bV-\bA\left(1-\frac{\delta_*}{|\bA|}\right)^+\right):\left(\bD-\bA\right)\geq0.$$

Choosing $\bA:= \bD\pm\varepsilon \bC$, with arbitrary $\bC\in L^2 (U)^{3\times 3}$ and $\varepsilon>0$, and after  the limit as $\varepsilon\to 0$ we obtain 

$$\int_U \bC: \left(\bV-\bD\left(1-\frac{\delta_*}{|\bD|}\right)^+\right)=0$$
for any  $\bC$, which implies \eqref{C3}.
\chiu

Note that here we provided a proof of \eqref{C1}, which is simplified and shorter than the one given in  \cite[Proposition 5.3]{old}.

\section{Approximations}\label{sectionApprox}
In this section, we prepare all the needed tools in order to prove Theorem \ref{main-thm1}. 
For any $n\in \mathbb{N}$, we introduce the following approximating system

\be\label{system22}
\ba{l}\displ
\diver\bv = 0 \textrm{ in } Q, \\ \displ
\pa_t\bv + \diver(\bv\otimes \bv) G_n(|\bv|^2) -\diver  \bS +\nabla p =  \bb \textrm{ in } Q, \\\displ
\pa_t \pf+ \bv \cdot \na \pf-  \Delta \pf=-\diver \bb +\bv\cdot\nabla p_s \textrm{ in } Q,
\ea
\ee

where  $G_n:\R\to \R$ is a  smooth function such that $G_n(u)=1$ if $|u|\leq n$, $G_n(u)=0$ if $|u|\geq 2n$ and $|G'_n|\leq \frac{2}{n}$. Next, we consider the following regularization of the constitutive equations (both in the bulk and on the boundary)

\begin{align}\label{Sn}
&\begin{array}{l}\displ\vspace{6pt}
\!\!\!\!\bS=\mathcal{S}_n(\pf, \bD {\bv})= \mathcal{Z}_n(\pf, \bD{\bv})+\left(1-\frac{\delta_*}{|\bD{\bv}|}\right)^+\bD{\bv} \ \ \ \mbox{ where }\mathcal{Z}_n(\pf, \bD{\bv}) := \tau(\pf)\frac{\bD{\bv}}{|\bD{\bv}|+\frac{1}{n}} \\\hfill \displ \vspace{6pt} \mbox{ with } \tau(\pf)=(p_s-\pf)^+ \textrm{ in } Q,
\end{array}\\ \label{sn}
&\!\!\bs= {s}_n(\bv_\tau)=\zeta_n(\bv_\tau) +\left(1-\frac{\beta_*}{|\bv_\tau|}\right)^+\bv_\tau \ \ \ \mbox{ where }
\zeta_n(\bv_\tau)= s_* \frac{\bv_{\tau}}{|\bv_{\tau}|+\frac{1}{n}} \mbox{ on } \Sigma,
\end{align}

and we complete the problem with boundary and  initial  conditions

\begin{align}\label{bc-n}
\bv\cdot\bn&=0, \na \pf \cdot \bn=0 \mbox{ on }  \Sigma, \\
\!\!\bv(0)&=\bv_0, \ \, \pf(0)={p}_0  \mbox{ in }  \O.\label{ic-n}
\end{align}

Note that both mappings $\bD\mapsto \mathcal{Z}_n(\pf, \bD)$ and $\bD\mapsto \left(1-\frac{\delta_*}{|\bD|}\right)^+\bD$ are monotone, i.e. 
  \begin{equation}\label{monSp1}(\bZ-\hat{\bZ}):(\bD-\hat{\bD})\geq 0 \mbox{ for any }\bZ=\mathcal{Z}_n(\pf, \bD), \hat{\bZ}=\mathcal{Z}_n(\pf, \hat{\bD}),
\end{equation}
 see formula (5.2) in \cite{old}, and  
\begin{equation}\label{monSp2}
(\bV-\hat{\bV}):(\bD-\hat{\bD})\geq 0\mbox{ for any } \bV=\left(1-\frac{\delta_*}{|\bD|}\right)^+\bD,  \hat{\bV}=\left(1-\frac{\delta_*}{|\hat{\bD}|}\right)^+\hat{\bD},
\end{equation}
see Lemma B.1 in \cite{BleMalRaj}. Therefore, due to the presence of the truncation in the convective term and the introduced approximations in the constitutive equations, the existence of weak solutions to system \eqref{system22}--\eqref{ic-n} can be proved through standard techniques of monotone operators, following also the spirit of the proof in \cite[Proposition 5.1]{old}. We enunciate the relevant result below and for the reader's convenience the proof can be found in  Appendix. 

\begin{prop}\label{prop}
Let $n\in \N$ be fixed. For any 
$$ \bv_0\in \Lnd{2}, \ p_0\in L^2(\O),\ \bb\in L^2(Q)  \mbox{ and } p_s\in L^5(Q),$$
 there exists a weak solution to the problem \eqref{system22}--\eqref{ic-n}. More precisely, for each $n\in\N$ there is a quadruplet $(\bv,\pf, \bS, \bs):= (\bv^n, \pf^n, \bS^n, \bs^n) $ such that
\begin{align}
&\bv \in  L^\infty (0,T; \Lnd{2})\cap L^2(0,T;W^{1,2}_{{\bn},{\rm div}}), \ \ \partial_t\bv\in {(L^2(0,T;W^{1,2}_{\bn}))}^*, \label{an1}\\
&\pf \in  L^\infty (0,T; L^2(\O))\cap L^2(0,T; W^{1,2}(\O)),\ \ \partial_t \pf \in (L^4(0,T; W^{1,2}(\O)))^*,\label{an2}\\
&\bS\in  L^{2}(Q)^{3\times 3}, \ \ \bs\in L^\frac{8}{3}(\Sigma)^3, \label{an3}
\end{align}

satisfying

\begin{align}
&\begin{array}{l}\displ\vspace{6pt}\!\!\!\! \!\!\!\int_0^T\!\!\! \!\langle \partial_t\bv, \bw\rangle +\!\! \int_{Q}\!\! (\bS\!:\! \bD\bw
 + G_n(|\bv|^2)\!\diver(\!\bv\otimes \bv\!)\! \cdot \!\bw) 
 +\!\! \int_{\Sigma}\!\! \bs\!\cdot\! \bw_{\btau} 
 = \!\!\int_{Q}\!\!\! \bb\cdot  \bw \textrm{ for all } \bw\!\in\! L^2(0, T; \!W^{1,2}_{\bn, \diver}),\end{array}
\label{an4}
 \\
&\label{an5}
\!\!\! \!\int_0^T\!\!\!   \langle \partial_t \pf, z\rangle  -\int_{Q}\!\!\pf \bv\!\cdot\! \nabla z  +\int_{Q}\!\! \nabla \pf\!\cdot\!\nabla z= \int_{Q}\!\!(\bb\cdot \nabla z- p_s\bv\!\cdot\!\nabla z)   \ \ 
\mbox{ for all } z\!\in\! L^4(0,T; W^{1,2}(\O)),
\end{align}
 
where
\begin{align}
& \bS=\mathcal{S}_n(\pf, \bD {\bv}) \textrm{ a.e. in } Q, \label{an6}\\
&\bs= s_n(\bv_\tau) \textrm{  a.e. in } \Sigma, \label{an7}
\end{align}

and 
\begin{equation}
\lim_{t\rightarrow 0^+}\|\bv(t)-\bv_0\|_2=0,  \ \lim_{t\rightarrow 0^+}\|\pf(t)-p{_0}\|_2=0. \label{an8}
\end{equation}
\end{prop}

\section{Proof of Theorem \ref{main-thm1}}\label{sectionProof}The proof is split in the following steps.\\
{\bf Step 1.}\textit{ Approximations.} From Proposition \ref{prop} and following the reconstruction of the pressure in \cite[Theorem 4.1, Step 2 of the proof]{old}, we get for each $n\in \N$ the existence of $(\bv^n,\pf^n, p^n, \mathbb{S}^n, \bs^n)$, with $p^n\in L^2(Q)$, satisfying

\begin{align}
&\ba{l}\displ\vspace{6pt}
\int_0^T\!\!\! \langle\partial_t \bv^n, \bw\rangle + \int_{Q}(\bS^n:\bD\bw+\diver(\bv^n\otimes \bv^n) G_n(|\bv^n|^2)\cdot \bw) +\int_{\Sigma}\!\!\bs^n\cdot \bw_{\btau} 
 = \int_{Q}p^n \diver \bw   \\ \displ \vspace{6pt} 
\hfill + \int_{Q} \!\!\!\bb\cdot \bw \ \ \textrm{ for all } \bw \in L^2(0, T; W^{1,2}_{\bn}),
 \label{weak-n}
 \ea\\
&\ba{l}\displ\vspace{6pt}
\int_0^T\!\!\langle \partial_t \pf^{n}, z\rangle - \int_{Q} \! (\pf^{n} \bv^{n}) \! \cdot\! \nabla z +\int_{Q}\! \nabla \pf^n\!\cdot\!\nabla z
=\int_{Q}\!(\bb \!-\!  p_s \bv^n)\! \cdot\! \nabla z\ \ 
\mbox{ for all } z \in L^4(0, T; W^{1,2}(\O)\!),
\ea\label{534}
\end{align}

with $\bS^n, \bs^n$ fulfilling \eqref{an6}, \eqref{an7} respectively, and satisfying also \eqref{an8}.  
\\

{\bf Step 2.} \textit{ Uniform estimates with respect to $n$ and limit as $n\to+\infty$.}
Setting $\bw:=\bv^n$  in  \eqref{weak-n} and $z:=\pf^n$ in  \eqref{534},  following the analogous step as in the proof of \cite[Theorem 4.1]{old}, we obtain 

\begin{align}\label{unif-n}
&\sup_{n}\left(\|\bv^n\|_{L^\infty(0,T; L^2(\O)^3)} + \|\bD \bv^n\|_{2, Q} \right) <+\infty,\\
\label{unif-pf}
&\sup_n \left(\|\pf^n\|_{L^\infty(0,T; L^2(\O))} + \|\nabla \pf^n \|_{2, Q}\right) <+\infty,\\
&\sup_n\left( \|\bv^n\|_{\frac{10}{3}, Q}+\|\pf^n\|_{\frac{10}{3}, Q}+\|\bv^n\|_{\frac{8}{3}, \Sigma}\right)<+\infty, \label{new0}\\
& \label{newe1}
\sup_n\left(\|\bZ^n\|_{\frac{10}{3}, Q}+\|\bV^n\|_{2, Q}+\|\bs^n\|_{\frac{8}{3},\Sigma}\right)<+\infty,
\end{align}

where we set $\bV^n:=\left(1-\frac{\delta_*}{|\bD{\bv^n}|}\right)^+\bD{\bv^n}$. Consequently, as $\sup_n \|G_n(|\bv^n|^2)\|_{L^\infty(Q)}\leq 1$ by employing \eqref{unif-n}, \eqref{new0} and Korn's inequality,  it follows that

\begin{equation} \label{newe2}
\sup_n \|G_n(|\bv^n|^2)\diver (\bv^n\otimes \bv^n)\|_{L^{\frac{5}{4}}(Q)}<+\infty.
\end{equation}

Now, let us introduce

 $$p_2^n:= (-\Delta_N )^{-1} \left(G^n(|\bv^n|^2)\diver (\bv^n \otimes \bv^n) \right),\ \ \ p_1^n:= p^n- p_2^n,$$
 
then

\begin{equation}\label{decp} \sup_n \left(\|p_2^n\|_{L^\frac{5}{4} (0,T; W^{1,\frac{5}{4}} (\O))}+ \|p_1^n\|_{L^2(Q)}\right) < +\infty,
\end{equation}

and this implies that

\be \label{unif-vt}
\sup_n\|\partial_t \bv^n\|_{{(L^2(0,T;W^{1,2}_{\bn})\cap L^5(Q)^3)}^*}<+\infty.
\ee

Analogously

\be\label{pf-time}
\sup_n\|\partial_t \pf^n\|_{{(L^4(0,T; W^{1,2}(\O)))}^*}<+\infty.
\ee

Then, there exist subsequences of $\{\bv^n\}$, $\{\pf^n\}$, $\{\bZ^n\}$, $\{\bV^n\}$, $\{\bs^n\}$, $\{p^n_1\}$, $\{p_2^n\}$, which we do not relabel, that  converge weakly and *-weakly in the corresponding function spaces. By virtue of the established limits, by the Aubin-Lions compactness lemma and the compact embedding of the Sobolev spaces into the space of traces, we also have
\begin{align}
&\pf^n  \to \pf  \textrm{ strongly in } L^q(Q) \mbox{ for all } q\in\left[1, \frac{10}{3}\right),\label{strongpf}\\
&\bv^n  \to \bv  \textrm{ strongly in } L^q(Q)^3 \mbox{ for all } q\in\left[1, \frac{10}{3}\right),\label{strongvi}\\
&\bv^n_\tau \to \bv_\tau \textrm{ strongly in } L^q(\Sigma)^3   \mbox{ for all } q\in\left[1, \frac{8}{3}\right).\label{vntau}
\end{align}
Since
$$\|G_n(|\bv^n|^2)\|_{L^\infty(Q)}\leq\! 1 \mbox{ and } G_n(|\bv^n|^2)\!\to\!\! 1 \mbox{ strongly in } L^q(Q) \mbox{ for all } q\!\in\! [1,+\infty),$$
it follows from \eqref{strongvi} that 
\begin{equation}
G_n(|\bv^n|^2)\diver (\bv^n\otimes \bv^n)\rightharpoonup \diver (\bv\otimes \bv) \textrm{ weakly in } L^{\frac{5}{4}}(Q)^3.\label{div-unif}
\end{equation}

Finally, with the obtained convergences it is standard to prove that
\be\label{a1}
\begin{split}\displ
\int_0^T \langle\partial_t \bv, \bw\rangle + \int_{Q}(\bZ+\bV): \bD\bw + \int_{Q}\bs\cdot\bw_{\btau}
 -\int_{Q}(\bv\otimes \bv): \bD\bw  = \int_{Q}p_1 \diver \bw\\ \displ
 -\int_{Q}\nabla p_2\cdot\bw +\int_{Q}\bb\cdot \bw\ \  \mbox{ for all }\bw\! \in\! L^2(0,T;W^{1,2}_{\bn})\cap L^5(Q)^3
 \end{split}
\ee
and
\be\label{}
\int_0^T\!\!\langle \partial_t \pf, z\rangle -  \int_{Q}(\pf \bv\cdot\nabla z  +  \int_{Q}\nabla \pf:\nabla z =\!\int_{Q}(\bb-p_s\bv)\!: \!\na z \ \ \mbox{ for all } z\!\in\! L^4(0,T; W^{1,2}(\O)).
\ee
 
{ \bf Step 3. }\textit{ Attainment of the constitutive equations on the boundary.} 
Using that 
$$\bs^n\rightharpoonup \bs \textrm{ weakly in } L^{\frac{8}{3}} (\Sigma)$$
and \eqref{vntau}, it easily follows
$$\limsup_{n\to + \infty}\int_{\Sigma} \bs^n\cdot\bv^n_\tau= \int_{\Sigma} \bs\cdot\bv_\tau.$$
Thus, a suitable adjustment of Proposition \ref{convergence-lemma} implies that \eqref{obiettivo-boundary} is fulfilled.
\\
 
{\bf Step 4. }\textit{Attainment of the constitutive equations in the bulk.}
In order to employ Proposition \ref{convergence-lemma} we need to prove the limsup property \eqref{limsupZ}, but as the solution itself can not be used as test function in \eqref{a1},  we follow the strategy as in \cite{old} and perform the $L^\infty$-truncation method. To this aim, we introduce
$$\bw^n:=T_{\lambda_n} (\bv^n-\bv):= (\bv^n-\bv) \min{\left\{1, \frac{\lambda^n}{|\bv^n-\bv|}\right\}}$$
where $\lambda^n \in [A,B]$ with $0<A< B<\infty$ will be suitably chosen numbers independent of $n$, but depending on parameter $N$ tending to $+\infty$, see details below. For the reader's convenience, we recall all the properties of $\bw^n$ below,
 \begin{align}\label{strongw}
&\bw^n\to 0 \ \mbox{strongly in} \ L^s(Q)\ \mbox{for every} \ s \in [1, +\infty),\\
\label{strongpw}
&\bw^n\to 0 \ \mbox{strongly in} \ L^2(\Sigma),\\
\label{weakw} &\bw^n \rightharpoonup 0\ \mbox{weakly in} \ L^2 (0,T; W^{1,2}_{\bn}),\\
\label{divwn}
&|{\rm div} \ \bw^n| \leq \begin{cases}\displ 0 & \mbox{if} \ |\bv^n - \bv|\leq\lambda^n\\ \displ
\frac{2 \lambda^n \left(|\nabla \bv^n| + |\nabla \bv|\right)}{|\bv^n - \bv|} & \mbox{if} \ |\bv^n-\bv|> \lambda^n,
\end{cases}\\
\label{nablawn}
& \nabla\bw^n\!\!=\!\!\begin{cases}\nabla \bv^n-\nabla \bv \ \ \ \ \   \mbox{if} \ |\bv^n - \bv|\leq\lambda^n\\ \displ
\begin{array}{l}\displ
\!\!\frac{\lambda^n}{|\bv^n - \bv|}(\nabla \bv^n-\nabla \bv)-\!\lambda^n(\bv^n-\bv)\!\otimes\! \frac{(\nabla \bv^n-\nabla \bv) (\bv^n-\bv)}{|\bv^n-\bv|^3} \\\displ \hfill
\mbox{if} \ |\bv^n-\bv|> \lambda^n.\end{array}
\end{cases}
\end{align}
Inserting $\bw^n$ as test function in \eqref{weak-n}, using the properties of  $\bw^n$
we get (cfr. \cite{old})
\begin{equation}\label{consw} 
\limsup_{n\to \infty}\int_{Q} (\bZ^n+\bV^n): \bD{\bw^n}\leq \limsup_{n\to\infty} \int_{Q} |p_1^n| |\diver \bw^n|.
\end{equation}
Now, let us define
\be\label{Zbar}\overline{\bZ}=
\begin{cases}\vspace{6pt}\displaystyle
\bO \ &\mbox{if} \  \bD\bv=\bO ,\\\vspace{6pt}\displaystyle
\tau(\pf)\frac{ \bD{\bv}}{|\bD{\bv}|}  &\mbox{if} \ \bD\bv\neq\bO. \end{cases}
\ee and
\begin{equation}
\label{Vbar}\overline{\bV}=\left(1-\frac{\delta_*}{|\bD{\bv}|}\right)^+\bD{\bv}.
\end{equation}
Employing  \eqref{weakw} and \eqref{divwn}, formula \eqref{consw} can be rewritten as
\begin{equation}
 \limsup_{n\to \infty}\!\! \int_{Q}\!\! (\bZ^n-\overline{\bZ})\!:\!\bD{\bw^n}+\!\int_{Q}\!\! (\bV^n-\overline{\bV})\!:\! \bD{\bw^n}\leq \limsup_{n\to\infty} \int_{\{|\bv^n-\bv|>\lambda^n\}} \!\!\!|p_1^n|\left(|\nabla \bv^n|\!+\! |\nabla \bv|\right) \frac{\lambda^n}{|\bv^n-\bv|}.
\end{equation}
Moving the part of the integral on the left-hand side on the set $\{|\bv^n-\bv|>\lambda^n\}$ to the right, it follows that
\begin{equation}\begin{split}\label{good1}
 \limsup_{n\to \infty} \int_{\{|\bv^n-\bv|\leq \lambda^n\}}\!\!\! \left({(\bZ^n-\overline{\bZ}): (\bD{\bv^n}-\bD{\bv})} + {(\bV^n-\overline{\bV}): (\bD{\bv^n}-\bD{\bv})}\right) \\ 
\leq C \limsup_{n\to \infty}\! \int_{\{|\bv^n-\bv|>\lambda^n\}} {I^n \frac{\lambda^n}{|\bv^n-\bv|}}\end{split}\end{equation}
where
 $$I^n:= |p_1^n|^2+  |\bZ^n|^2 + |\overline{\bZ}|^2 + |\bV^n|^2+|\overline{\bV}|^2+ |\nabla \bv^n|^2 + |\nabla \bv|^2.$$
 
Note that it holds (see formula (6.60) in \cite{old})
\begin{equation}\label{good2a}
\left((\bZ^n-\overline{\bZ}):(\bD{\bv^n}- \bD{\bv})\right)^- \to 0\ \mbox{strongly in} \ L^1(Q),
\end{equation}
and analogously the monotonicity implies 
\begin{equation}
\label{good2b}\left((\bV^n-\overline{\bV}):(\bD{\bv^n}- \bD{\bv})\right)^- \to 0\ \mbox{strongly in} \ L^1(Q).\end{equation}
Let $N\in \N$ and fix $A=N$, $B=N^{N+1}$. For $i\in \{1,...,N\}$ let us define
$$Q_i^n:= \{N^i\leq|\bv^n-\bv|\leq N^{i+1}\}.$$
Since $$\sum_{i=1}^N \int_{Q_i^n} I^n \leq C^*,$$
for every $n$ there exists $i_n\in \{1,...,N\}$ such that $$\int_{Q_{i_n}^n} I^n \leq \frac{C^*}{N}.$$ 
Set $\lambda^n:= N^{i_n}$, then it holds
 $$ \int_{\{|\bv^n-\bv|> \lambda^n\}} I^n \frac{\lambda^n}{|\bv^n-\bv|}=  \int_{Q_{i_n}^n} I^n \frac{N^{i_n}}{|\bv^n-\bv|} +  \int_{\{|\bv^n-\bv|> N^{i_n+1}\}} I^n \frac{N^{i_n}}{|\bv^n-\bv|} \leq \frac{C^*}{N}$$
 where we keep the symbol $C^*$ for a different constant.
The latter relation,  \eqref{good1}, \eqref{good2a} and \eqref{good2b} give
\begin{equation}\label{fin1} 
\begin{split}\vspace{6pt}\displ
\limsup_{n\to +\infty} \left(\int_{\{|\bv^n-\bv|\leq \lambda^n\}}\!\!\!\!\! |(\bZ^n-\overline{\bZ}):(\bD{\bv^n}-\bD{\bv})|\right. \\ \vspace{6pt}\displ
\left.+ \int_{\{|\bv^n-\bv|\leq \lambda^n\}}\!\!\!\! |(\bV^n-\overline{\bV}):(\bD{\bv^n}-\bD{\bv})| \right)
\leq \frac{C^*}{N}.
\end{split}\end{equation} 
Using that
\begin{equation}\label{fin4}
\begin{split}
\int_{Q} \sqrt{|(\bZ^n-\overline{\bZ}):(\bD{\bv^n}-\bD{\bv})|} = 
\int_{\{|\bv^n-\bv|\leq N\}} \sqrt{|(\bZ^n-\overline{\bZ}):(\bD{\bv^n}-\bD{\bv})|} \\ 
+ \int_{\{|\bv^n-\bv|> N\}} \sqrt{|(\bZ^n-\overline{\bZ}):(\bD{\bv^n}-\bD{\bv})|}
\end{split}\end{equation}
and
\begin{equation}\label{fin5}\begin{split}
\int_{
Q} \sqrt{|(\bV^n-\overline{\bV}):(\bD{\bv^n}-\bD{\bv})|} = 
\int_{\{|\bv^n-\bv|\leq N\}} \sqrt{|(\bV^n-\overline{\bV}):(\bD{\bv^n}-\bD{\bv})|}  \\
+ \int_{\{|\bv^n-\bv|> N\}} \sqrt{|(\bV^n-\overline{\bV}):(\bD{\bv^n}-\bD{\bv})|},
\end{split}\end{equation}
by H\"{o}lder's and Chebyshev's inequalities we obtain
\begin{equation}
\limsup_{n\to +\infty}\int_{Q} \sqrt{|(\bZ^n-\overline{\bZ}):(\bD{\bv^n}-\bD{\bv})|}  \leq \frac{2\overline{C}}{\sqrt{N}},
\end{equation}
and 
\begin{equation}
\limsup_{n\to +\infty}\int_{
Q} \sqrt{|(\bV^n-\overline{\bV}):(\bD{\bv^n}-\bD{\bv})|}\leq \frac{2\overline{C}}{\sqrt{N}},
\end{equation}
which means, by letting $N\to\infty$, that
$$(\bZ^n-\overline{\bZ}):(\bD{\bv^n}-\bD{\bv})\to 0 \ \mbox{a.e. in} \ Q,$$
$$(\bV^n-\overline{\bV}):(\bD{\bv^n}-\bD{\bv})\to 0 \ \mbox{a.e. in} \ Q.$$
Egoroff's theorem then gives that for all $\varepsilon>0$ there exists $U\subset Q$, $|Q\setminus U|\leq \varepsilon$ such that
\begin{equation}\label{fin6} \int_U {(\bZ^n-\overline{\bZ}):(\bD{\bv^n}-\bD{\bv})}\to 0 , 
\end{equation}
\begin{equation}\label{fin7}\int_U {(\bV^n-\overline{\bV}):(\bD{\bv^n}-\bD{\bv})}\to 0.  \end{equation}
We conclude, thanks to the weak convergences of $\bZ^n, \bV^n, \bD{\bv^n}$ respectively to $\bZ, \bV, \bD{\bv}$
 that
$$\lim_{n\to\infty}\int_U \bZ^n: \bD{\bv^n} = \lim_{n\to \infty} \int_U \bZ^n: \bD{\bv}= \int_U \bZ: \bD{\bv}$$
and
$$\lim_{n\to\infty}\int_U \bV^n: \bD{\bv^n} = \lim_{n\to \infty} \int_U \bV^n: \bD{\bv}= \int_U \bV: \bD{\bv}.$$
Finally, all assumptions of  Convergence Lemma \ref{prop} are fulfilled, thus \eqref{eqS} holds a.e. in $U$. Since $|Q\setminus U|\leq \varepsilon$ we can let $\varepsilon\to 0$ and obtain that \eqref{eqS} holds a.e. in $Q$.
Theorem \ref{main-thm1} is proved.

\section{Appendix}
Our goal is to prove Proposition \ref{prop}. Let us recall that in this section we fix $n\in \mathbb{N}$, we consider $G_n$ smooth function with the properties stated at the beginning of Section \ref{sectionApprox} and the regularization of the material responses
given in \eqref{Sn} and \eqref{sn}. In what follows, to simplify the notation we drop the indices $n$. 

\subsection*{\bf Proof of Proposition \ref{prop}}
The proof is split in the following steps.
\\

{\bf Step 1.}{\textit{ Approximations.}
For any $m\in \mathbb{N}$, we look for 
 \be
 \bv^{m}(t,x):=\sum_{r=1}^m c_r^{m}(t) \bw^r(x), \ \ \ \pf^{m}(t, x):=\sum_{r=1}^m d_r^{m}(t) z^r(x)
\ee
satisfying
\begin{equation}
\label{alpha}\ba{l}\displ \vspace{6pt}
\left(\! \frac{d \bv^{m}}{dt}, \bw^r\!\!\right)\!+\! (\mathbb{S}^m, \bD \bw^r ) 
 +( \diver (\bv^{m}\!\otimes\!\bv^{m})G(|\bv^{m}|^2),\bw^r) \\\displ \vspace{6pt}
 +  ( \bs^{m}, \bw^r)_{\pa\O} = (\bb,\bw^r) 
 \ \ r=1,\dots, m
\ea
\end{equation}
where
\begin{align}
&\begin{array}{l}\displ\vspace{6pt}
\!\!\mathbb{S}^m:=\mathcal{S}(\pf^m, \bD\bv^m)\!=\!\tau(\pf^m)\frac{\bD{\bv^m}}{|\bD{\bv^m}|+\frac{1}{n}}+\bD{\bv^m}\!\!\left(\!\!1-\!\frac{\delta_*}{|\bD{\bv^m}|}\!\right)^+ \\ \displ\vspace{6pt}\hfill
\mbox{ with } \tau(\pf^m)\!=\!(p_s-\pf^m)^+,
\end{array}\\
&\bs^m:=s(\bv_\tau^m)\!=s_* \frac{\bv_{\tau}^m}{|\bv_{\tau}^m|+\frac{1}{n}} +\bv_\tau^m\left(1-\frac{\beta_*}{|\bv_\tau^m|}\right)^+,
\end{align}
and 
\begin{equation}\label{beta}
\ba{l}\displ  
\left( \partial_t \pf^{m}, z^r \right) \!-\!  ( \pf^{m}\bv^{m},\nabla z^r ) \!+\!  ( \nabla\pf^{m},  \nabla z^r)  
= (\bb, \nabla z^r)
 -(p_s\bv^{m},\nabla z^r ) 
\\\displ 
\hfill  r=1,\dots, m,
\ea\end{equation}
where $\{\bw^i\}_{i\in\N}$ is an orthogonal basis in $W_{\mathbf{n}, \diver}^{1,2}$ consisting of eigenfunctions of the Stokes operator with boundary conditions $\bw^i\!\cdot\! \bn\!=\!0$ and $[(\bD\bw^i) \bn]_\tau=\!\bnul$ on $\partial\Omega$, while  $\{z_j\}_{j\in \mathbb{N}}$ is an orthogonal basis in $W^{1,2}(\Omega)$ consisting of  eigenfunctions of the Laplace operator subject to the Neumann homogeneous boundary conditions. The system is supplemented with the corresponding initial conditions $\bv^m_0$ and $p_0^m$, obtained by projection $ \bv_0\in \Lnd{2}$ onto the span of $[ \bw^1,\dots ,\bw^m]$ and respectively $p_0\in L^2(\O)$ onto the span of $ [z^1,\dots, z^m]$. Then the local in time existence of $\bv^{m}$ and $\pf^{m}$ follows from the Caratheodory theory for systems of ordinary differential equations, whereas the global in time existence is a consequence of the uniform estimates established below. 

{\bf Step 2. }\textit{Uniform estimates.} Multiplying \eqref{alpha} by $c^m_r(t)$ and \eqref{beta} by $d^m_r(t)$ and taking the sum over $r=1, \dots, m$, we obtain
\begin{align}
\label{gamma}
&\ba{l}\displ\vspace{6pt}
\frac{1}{2} \frac{d}{dt}\|\bv^{m}(t)\|_2^2   +  \int_{\{|\bD \bv^m|>\delta_*\}}|\bD \bv^m|^2
+ \int_{\Omega}\!\!  \tau(\pf^m)\frac{|\bD{\bv^m}|^2}{|\bD{\bv^m}|+\frac{1}{n}}
\\  \displ\vspace{6pt}\hfill+\int_{\partial \Omega}\!\!\! \bs^{m}\! \cdot\! \bv_{ \tau}^{m}=( \bb, \bv^{m})+ \int_{\{|\bD \bv^m|>\delta_*\}} \delta_*|\bD \bv^m|,
\ea
\\
\label{delta2}
&\frac{1}{2} \frac{d}{d t} \|\pf^{m}(t)\|_2^2 +  \| \nabla p_{f}^{m}(t)\|_2^2 = (\bb, \nabla\pf^m)-(p_s\bv^m, \nabla \pf^m),
\end{align}
Adding $\int_{\{|\bD \bv^m|\leq\delta_*\}}|\bD \bv^m|^2$ to both sides of \eqref{gamma} 
\begin{equation*}
\frac{1}{2} \frac{d}{dt}\|\bv^{m}(t)\|_2^2   +  \int_{\O}|\bD \bv^m|^2
+ \int_{\Omega}\!\! \tau(\pf^m)\frac{|\bD{\bv^m}|^2}{|\bD{\bv^m}|+\frac{1}{n}}
+\int_{\partial \Omega}\!\!\! s_n(\bv_\tau^{m})\! \cdot\! \bv_{ \tau}^{m}\leq
(\bb, \bv^{m}) + \int_{\Omega}\delta_*|\bD \bv^m|+ \delta_*^2|\Omega|
\end{equation*}
and then by Young's inequality, we get  
\begin{equation}\label{gamma2}
\frac{1}{2} \frac{d}{dt}\|\bv^{m}(t)\|_2^2   + \frac{1}{2} \int_{\O}|\bD \bv^m|^2
+ \int_{\Omega}\!\! \tau(\pf^m)\frac{|\bD{\bv^m}|^2}{|\bD{\bv^m}|+\frac{1}{n}}
+\int_{\partial \Omega}\!\!\! s_n(\bv_\tau^{m})\! \cdot\! \bv_{ \tau}^{m}\leq
(\bb, \bv^{m}) +\frac{3}{2} \delta_*^2|\Omega|.
\end{equation}
Integrating in time, by Korn's and Young's inequalities, using also the fact that the last two terms on the left-hand side of \eqref{gamma2} are non-negative, one concludes that
\begin{equation}\label{AEO}
\sup_{t\in [0,T]} \|\bv^{m} (t)\|_2^2 +  \int_{Q} |\bD {\bv^{m}}|^2
\leq C(\bb, \bv_0, \delta_*, |Q|).
\end{equation}
By the interpolation inequality
\begin{align}
&\|z\|_{\frac{10}{3}}\leq \|z\|_2^{\frac{2}{5}}\|z\|_6^{\frac{3}{5}}\leq C\|z\|_2^{\frac{2}{5}}\|z\|_{1,2}^{\frac{3}{5}}\label{dt0}
\end{align}
and the trace inequalities (see \cite[Lemma 1.11]{BMR}), we obtain
\begin{equation}\label{AE2}
\sup_m \left(\!\|\bv^m\|_{\frac{10}{3}, 
Q} + \|\bv^m_\tau\|_{\frac{8}{3}, \Sigma}\!\!\right)<+\infty.
\end{equation}
As a consequence integrating in time \eqref{delta2} we deduce that
\begin{equation}\label{AE1}
\sup_{t\in [0,T]}\|\pf^{m}(t)\|_2^2+ \int_{Q}\!\!|\nabla \pf^{m}|^2 \leq C\|\bb\|_{2,Q}^2+ C\|p_s\|_{5,Q}^2 \|\bv^m\|_{\frac{10}{3},Q}^2+ \|{p}_0 \|_2^2
\end{equation}
and thus
\begin{equation}\label{AE3}
\sup_{t\in(0,T)}\|\pf^m (t)\|_2 + \|\nabla \pf^m\|_{L^2(Q)}\!\!\leq C(\bb, p_s, p_0).
\end{equation}
Again \eqref{dt0} gives
\begin{equation}\label{p103} \sup_m \|\pf\|_{\frac{10}{3}, Q} < +\infty.
\end{equation}
Recalling the explicit formulas for $\bS^m$ and $\bs^m$ it then follows 
\begin{equation}\label{AE4}
\sup_m\left(\|\bS^m\|_{2, Q} +\|\bs^m\|_{\frac{8}{3}, \Sigma}\right)<+\infty.
\end{equation}
Employing the inequality
$$ \|z\|_{4}\leq \|z\|_2^{\frac{1}{4}}\|z\|_6^{\frac{3}{4}}\leq C\|z\|_2^{\frac{1}{4}}\|z\|_{1,2}^{\frac{3}{4}}$$
we deduce corresponding uniform estimates for $\bv^m$ and $\pf^m$ respectively in $L^4(Q)^3$ and $L^4(Q)$, then by virtue of them
it results
\begin{equation}\label{AE6}
\sup_m \|\partial_t\pf^m\|_{(L^4(0,T; W^{1,2}(\O)))^*}<+\infty.
\end{equation}
Analogously and by virtue of the truncation in the convective term, we also get 
\begin{equation}\label{AE5}
\sup_m \|\partial_t\bv^m\|_{{(L^2(0,T; W^{1,2}_{\bn, \diver}))}^*}<+\infty.
\end{equation}

{\bf Step 3. }\textit{Limit.} By virtue of uniform estimates established above there exist subsequences of $\{\bv^m\}, \{\pf^m\}$, $\{\bS^m\}$ and $\{\!\bs^m\!\}$, converging  respectively  weakly (or *-weakly) to $\bv, \pf, \bS$ and $\bs$ in the corresponding function spaces. Furthermore, Aubin-Lions compactness Lemma and its variant including the trace theorem imply the following strong convergences:
\begin{align}
&\bv^m\to \bv \mbox{ a.e. in } Q \mbox{ and strongly in } L^q(Q)^3 \mbox{ for any } q\in \left[1, \frac{10}{3}\right),\label{L1}\\
&\pf^m\to \pf  \mbox{ a.e. in } Q \mbox{ and strongly in } L^q(Q)^3 \mbox{ for any } q\in \left[1, \frac{10}{3}\right),\label{L2}\\
& \bv^m_\tau \to \bv_\tau  \mbox{ a.e. in } \Sigma \mbox{ and strongly in } L^q(\Sigma)^3 \mbox{ for any } q\in \left[1, \frac{8}{3}\right).\label{L3}
\end{align}
As a consequence $\bv, \pf, \bS$ and $\bs$ fulfill the weak formulations stated in Proposition \ref{prop}.
\\

{\bf Step 4. } \textit{Attainment of the constitutive equations.} 
The convergence
$$\bs^m \rightharpoonup \bs \mbox{ weakly in } L^{\frac{8}{3}}(\Sigma)$$
together with \eqref{L3} ensures that
\be
\lim_{m\to +\infty} \int_{\Sigma}\bs^m\cdot\bv_\tau^m=\int_{\Sigma}\bs\cdot\bv_\tau.
\ee
Then, thanks to the monotonicity it is standard to prove that 
$$\bs=s(\bv_\tau).$$
Next,  it follows from the monotonicity that
\begin{equation}\label{M1}
0\leq \!\int_{Q}\left(\bS^m- \mathcal{S}(\pf^m, \bA\!)\right)\!:\!\left(\bD\bv^m-\bA\right)\ \ \mbox{   for all } \bA \in L^2(Q).
\end{equation}
Now, note that by \eqref{L2},
\begin{equation}\label{abbatiel} \ba{c}\displ \vspace{6pt}
\mathcal{S}(\pf^m, \bA):=(p_s-\pf^m)^+\!\frac{\bA}{|\bA|+\frac{1}{n}}+\bA\left(1-\frac{\delta_*}{|\bA|}\right)^+\\ \displ\vspace{6pt}
\to(p_s-\pf)^+\!\frac{\bA}{|\bA|+\frac{1}{n}}+\bA\left(1-\frac{\delta_*}{|\bA|}\right)^+=:\mathcal{S}(\pf, \bA\!) \mbox{ strongly in } L^2(Q) \ea\end{equation}
while, as $\bv$ can play the role of a test function in the established weak formulation, it is standard to obtain
\begin{equation}\label{EI3}
\limsup_{m\to+\infty} \int_{Q}\!\!\bS^m\!\!:\!\bD\bv^m\leq \int_{Q}\!\!\bS\!:\!\bD\bv.
\end{equation}
 Finally, thanks to the convergences 
\begin{align*}
&\bD\bv^m\rightharpoonup\bD\bv \mbox{ weakly in } L^2(Q), \\
&\bS^m\rightharpoonup\bS \mbox{ weakly in } L^2(Q),
\end{align*}
and \eqref{abbatiel},
the limit as $m\to + \infty$ in \eqref{M1} gives
\begin{equation}\label{M2}
0\leq \!\int_{Q}\!\!\left(\bS- \mathcal{S}(\pf, \bA\!)\right)\!:\!\left(\bD\bv-\bA\right)\ \ \mbox{   for all } \bA \in L^2(Q).
\end{equation}
At this point, it is standard to choose  $\bA=\bD\bv\pm\varepsilon\bB$ for arbitrary $\bB\in L^2(Q)$  and $\varepsilon>0$ and arrive at
$$ 0= \!\int_{Q}\bB:\left(\bS- \mathcal{S}(\pf,\bD\bv)\right)\ \ \mbox{   for all } \bB \in L^2(Q),$$
which implies $\bS=\mathcal{S}(\pf,\bD\bv)$ a.e. in $Q$. The proof of Proposition \ref{prop} is complete.
\qed

\section*{Acknowledgments}
The research of A.~Abbatiello is supported by Einstein Foundation, Berlin.  A.~Abbatiello is also member of the Italian National Group for the Mathematical Physics GNFM-INdAM. M.~Bul\'{i}\v{c}ek, T.~Los, J.~M\'alek,  and O.~Sou\v{c}ek acknowledge support of the project 18-12719S financed by the Czech Science Foundation. T.~Los is also thankful to the institutional support  through the project GAUK 550218 and the Charles University Research Program No. UNCE/SCI/023.
The authors are also grateful to Hausdorff Research Institute for Mathematics in Bonn (Germany), where the work was partially performed. M.~Bul\'{i}\v{c}ek, J.~M\'alek,  and O.~Sou\v{c}ek are members of the Ne\v{c}as Center for Mathematical Modeling.


\begin{thebibliography}{1}

\bibitem{AbbFeir}
Abbatiello,~ A., Feireisl,~E.: \emph{On a class of generalized solutions to equations describing incompressible viscous fluids.}
Ann. Mat. Pura Appl.  (2019)  https://doi.org/10.1007/s10231-019-00917-x

\bibitem{old}
Abbatiello,~ A., Los,~T., M{\'a}lek, J., Sou\v{c}ek,~O.: \emph{Three--dimensional flows of pore pressure--activated {B}ingham fluids.} Math. Models Methods Appl. Sci. \textbf{29}, 2089--2125 (2019)

\bibitem{BleMalRaj}
Blechta,~J., M{\' a}lek,~J., Rajagopal,~K.R.: \emph{On the classification of
  incompressible fluids and a mathematical analysis of the equations that
  govern their motion.} To appear in SIAM J.  Math. Anal. Arxive Preprint Series \textbf{{\bf arXiv
  1902.04853v1}} (2019)

\bibitem{sol-lip}
Breit,~D., Diening,~L., S.~Schwarzacher,~S.: \emph{Solenoidal Lipschitz
  truncation for parabolic PDE's.} Math. Models Methods Appl. Sci. \textbf{23},
  2671--2700  (2013)

\bibitem{BM}
Bul{\'\i}{\v{c}}ek,~M., M\'alek,~J.: \emph{On unsteady internal flows of
  {B}ingham fluids subject to threshold slip on the impermeable boundary.}
  Recent Developments of Mathematical Fluid Mechanics, (Eds. H. Amann, Y. Giga,
  H. Okamoto, H. Kozono, M. Yamazaki) Birkh\"{a}user-Verlag (2014)

\bibitem{BMR}
Bul{\'\i}{\v{c}}ek,~M., M\'alek,~J., Rajagopal,~K.R.: \emph{Navier's slip and
  evolutionary {N}avier-{S}tokes-like systems with pressure and shear-rate
  dependent viscosity.} Indiana Univ. Math. J. \textbf{56}  no.~1,
  51--86 (2007)

\bibitem{CM}
Chupin,~L., Math\'{e},~J.: \emph{Existence theorem for homogeneous incompressible
  {N}avier-{S}tokes equation with variable rheology.} European Journal of
  Mechanics. B. Fluids \textbf{61}  no.~part 1, 135--143 (2017)

\bibitem{MaRuSh}
M\'alek,~J., R\r{u}\v{z}i\v{c}ka,~M., Shelukhin,~V.V.:
  \emph{Herschel-{B}ulkley fluids: existence and regularity of steady flows.}
  Math. Models Methods Appl. Sci. \textbf{15} no.~12, 1845--1861 (2005)
  
  \bibitem{RA} Rajagopal,~K.R.: \emph{On implicit constitutive theories.}
   Appl. Math. {\bf 48}  no.~4, 279--319 (2003)
   
   \bibitem{RA2} Rajagopal,~K.R.: \emph{On implicit constitutive theories for fluids.}
  J. Fluid Mech. {\bf 550}  243--249 (2006)
  
\bibitem{She}
Shelukhin,~V.V.: \emph{Bingham viscoplastic as a limit of non-{N}ewtonian
  fluids.} J. Math. Fluid Mech. \textbf{4}  no.~2, 109--127 (2002)
\end{thebibliography}
\end{document}